# ON THE STOCHASTIC CALCULUS METHOD FOR SPINS SYSTEMS

By Samy Tindel

*Université Henri Poincaré (Nancy)*

In this note we show how to generalize the stochastic calculus method introduced by Comets and Neveu [*Comm. Math. Phys.* **166** (1995) 549–564] for two models of spin glasses, namely, the SK model with external field and the perceptron model. This method allows to derive quite easily some fluctuation results for the free energy in those two cases.

**1. Introduction.** For $N \geq 1$, let $\Sigma_N = \{-1; 1\}^N$, and denote by $\sigma = (\sigma_1, \ldots, \sigma_N)$ a typical element of $\Sigma_N$. In this paper we will first consider the usual Sherrington–Kirkpatrick (SK) model with external field based on this space of configurations, that is, a random measure on $\Sigma_N$ induced by the following Hamiltonian:

$$(1) \qquad -H_N(\sigma) = \frac{\beta}{N^{1/2}} \sum_{i<j} g_{i,j} \sigma_i \sigma_j + h \sum_{i \leq N} \sigma_i.$$

In the previous equation, $\beta$ is a strictly positive parameter that stands for the inverse of the temperature of the system, $\{g_{i,j}; 1 \leq i < j\}$ is a family of i.i.d. standard Gaussian random variables, and $h$ is a strictly positive coefficient representing the external field, under which the spins tend to take the value $+1$. The measure under consideration on $\Sigma_N$ is then the measure $G_N$ whose density with respect to the counting measure is given by $Z_N^{-1} e^{-H_N(\sigma)}$, where $Z_N$ is the normalization constant given by

$$Z_N = \sum_{\sigma \in \Sigma_N} e^{-H_N(\sigma)}.$$

For $n \geq 1$ and a function $f : \Sigma_N^n \to \mathbb{R}$, we will set $\rho(f)$ for the average of $f$ under the measure $G_N^{\otimes n}$, that is,

$$\rho(f) = Z_N^{-n} \sum_{\sigma^1, \ldots, \sigma^n \in \Sigma_N} f(\sigma^1, \ldots, \sigma^n) \exp\left(-\sum_{l=1}^{n} H_N(\sigma^l)\right).$$

At high temperature (i.e., when $\beta$ is small enough), it has been shown (see, e.g., [8]) that this classical SK model could be understood in its essential features if one could study the asymptotic behavior of the random variable $Z_N$ and of the quantity (called overlap between $\sigma^1$ and $\sigma^2$)

$$R_{1,2} = \frac{1}{N} \sum_{i \leq N} \sigma_i^1 \sigma_i^2, \tag{2}$$

where $\sigma^1$ and $\sigma^2$ are taken as two independent configurations under $G_N$. Then, after the introduction and formalization of the cavity and smart path methods, many limit theorems have been obtained for $Z_N$ and $R_{1,2}$: we refer to [9] for the self averaging property for $R_{1,2}$ and the limit of $\frac{1}{N}\log(Z_N)$, to [5] (resp. [13]) for a central limit theorem for the fluctuations of $Z_N$ (resp. $R_{1,2}$), and to [10, 14] for a (quenched) large deviation principle for the overlap $R_{1,2}$.

On the other hand, it has been shown in [3] (and used extensively, e.g., in [2]) that in the case of the SK model without external field, that is, when $h = 0$ in (1), one could use some simple stochastic calculus tools in order to simplify the long calculations involved in the asymptotic results mentioned above. This is achieved by replacing the Gaussian path $t^{1/2}g$, $t \in [0,1]$, by a Brownian motion $B_t$, and applying then Itô's formula, that gives immediately the differential along this path. However, in case of the SK model with external field, the Gaussian paths are of the form

$$t^{1/2}g + r(1-t)^{1/2}X, \qquad t \in [0,1], \tag{3}$$

where $g$ and $X$ are two independent standard Gaussian random variables, and $r$ is a positive coefficient. The stochastic calculus analogous of this path would then be

$$B_t + r\beta_{1-t}, \qquad t \in [0,1],$$

where $B$ and $\beta$ are now two independent Brownian motions. Dealing with these two Brownian motions running in opposite directions has been seen as an obstacle to the generalization of the Comets–Neveu method to the case where $h \neq 0$. In this paper our aim is to show how to bypass this difficulty by just invoking the fact that $\{\beta_{1-t}; t \in [0,1]\}$ can also be seen as the solution to a stochastic differential equation on which we can perform an integration by parts. Then, we will show that the CLT for $Z_N$ can be obtained by applying Itô's formula to the fluctuations of this last quantity. This gives an alternative (and maybe easier) way to [5] to prove this result. Notice, however, that we do not try to use here the nice path introduced in [5], and, hence, the temperature region where our results hold is smaller than in this last reference. The reason why we do not do it is twofold: first, we wish to insist here on the simplicity of the tools we use rather than on the optimality of the result, and on the other hand, one of our aims is to generalize our method to other (and more complex) spin systems.

We will illustrate this second point by considering an analogous fluctuation result for the perceptron model. This spin glass system, motivated by some neural computing considerations, is still based on the configuration space $\Sigma_N = \{-1; 1\}^N$, and is induced by the random Hamiltonian

$$-H_{N,M}(\sigma) = \sum_{k \leq M} u\left(\frac{1}{N^{1/2}} \sum_{i \leq N} g_{i,k}\sigma_i\right),$$

where $M = \alpha N$ for a (small enough) proportional coefficient $\alpha$, $u$ is a bounded function and $\{g_{i,k}; i, k \geq 1\}$ is a family of independent standard Gaussian random variables. We refer to [6] for the computational motivation of this model, and to [13], Chapter 3, and [11] for the basic asymptotic results for the Gibbs measure induced by $H_{N,M}$. Since the limit theorems obtained for the perceptron model are based again on the analysis of the overlap quantity $R_{1,2}$ defined by (2), and on the use of some elaborated versions of the path (3), our stochastic calculus method still applies here, and, indeed, it will allow us to get quite easily a central limit theorem for the normalizing constant

$$Z_{N,M} = \sum_{\sigma \in \Sigma_N} e^{-H_{N,M}(\sigma)}.$$

Our paper is divided as follows: in the next section we will treat the Sherrington–Kirkpatrick case, for which we define the stochastic path in detail at Section 2.1, and then get the fluctuation result for $Z_N$ at Section 2.2. Section 3 is devoted to the perceptron case: at Section 3.1 we introduce the model and state our central limit theorem. Then, at Section 3.2, we apply the stochastic path method to prove this last result.

Notice that, throughout thee paper, $\kappa$ will designate a positive constant and $\mathcal{R}_N$, $\hat{\mathcal{R}}_N$, $\mathcal{R}_{m,N}$, and so on, some small remainders. The exact value of those quantities are generally irrelevant, and can change from line to line.

**2. The Sherrington–Kirkpatrick case.** In this section we will try to explain the stochastic calculus method and apply it to recover the fluctuations' results proved in [5] in a simple way.

2.1. *Definition of the path.* This section is devoted to the definition and some basic properties of the interpolating path we will consider throughout the paper: let $\{B_{i,j}; 1 \leq i < j\}$ and $\{W_i; i \geq 1\}$ be two collections of independent Brownian motions, and $\{\eta_i; i \geq 1\}$ a family of independent standard Gaussian random variables. All those objects will be defined on a complete stochastic basis $(\Omega, \mathcal{F}, (\mathcal{F}_t)_{t \in [0,1]}, P)$, and we will assume that all the $B_{i,j}, W_i$ are $\mathcal{F}_t$-adapted, and that all the random variables $\eta_i$ are $\mathcal{F}_0$-measurable. For $i \geq 1$, let $X_i$ be the unique solution (see, e.g., [7]) to the stochastic differential equation

$$(4) \qquad X_i(t) = \eta_i - \int_0^t \frac{X_i(s)}{1-s} \, ds + W_i(t), \qquad t \in [0,1].$$

It is easily checked that, setting $\hat{X}_i(t) = X_i(1-t)$, the process $\{\hat{X}_i(t); t \in [0,1]\}$ is a Brownian motion with respect to its natural filtration. Hence, $X_i$ can be seen as a reversed time Brownian motion. In the sequel, if $Y$ is a square integrable continuous semi-martingale, $\langle Y \rangle$ will stand for the quadratic variation process of $Y$. We will first label the following lemma for further use.

LEMMA 2.1. *Let $k \geq 1$, $X = (X_1, \ldots, X_k)$, where $X_i$ is the solution to (4), and $\varphi : \mathbb{R}^k \to \mathbb{R}$ be a $C^2$ function having at most exponential growth together with its first two derivatives. Set also $\eta = (\eta_1, \ldots, \eta_k)$. Then, for any $t \in [0,1]$,*

$$\mathbf{E}[\varphi(X(t))] = \mathbf{E}[\varphi(\eta)] - \tfrac{1}{2} \int_0^t \mathbf{E}[\Delta \varphi(X(s))] \, ds.$$

PROOF. According to Itô's formula, we have

$$(5) \qquad \varphi(X(t)) = \varphi(\eta) - \sum_{i=1}^k \int_0^t \partial_{x_i} \varphi(X(s)) \frac{X_i(s)}{1-s} \, ds + \sum_{i=1}^k \int_0^t \partial_{x_i} \varphi(X(s)) \, dW_i(s) + \frac{1}{2} \int_0^t \Delta \varphi(X(s)) \, ds.$$

Recall now that for a Gaussian vector $Y = (Y_1, \ldots, Y_k)$ and a $C^1$ function $\psi$ on $\mathbb{R}^k$ having at most exponential growth, we have, for $i \leq k$,

$$(6) \qquad \mathbf{E}[\psi(Y) Y_i] = \sum_{j=1}^k \mathbf{E}[Y_i Y_j] \mathbf{E}[\partial_{x_j} \psi(Y)].$$

In particular, since $X_i(s)$ is a $\mathcal{N}(0, 1-s)$ random variable independent from the other $X_j$'s, we have, for $s \in [0,1]$,

$$\mathbf{E}\left[\partial_{x_i} \varphi(X(s)) \frac{X_i(s)}{1-s}\right] = \mathbf{E}[\partial_{x_i^2}^2 \varphi(X(s))].$$

Taking expectations in (5) and applying this last identity, we get the desired result. □

Going back to our interpolating Hamiltonian, define, for $t \in [0,1]$,

$$(7) \quad -H_{N,t}(\sigma) = \frac{\beta}{N^{1/2}} \sum_{i<j} B_{i,j}(t)\sigma_i\sigma_j + \beta q^{1/2} \sum_{i \leq N} X_i(t)\sigma_i + h \sum_{i \leq N} \sigma_i,$$

where $q$ is the usual $L^2$-limit of the overlap, that is, the solution to the implicit equation

$$q = \mathbf{E}[\tanh^2(\beta q^{1/2}Y + h)],$$

where $Y$ stands for a standard Gaussian random variable. Let us write then

$$Z_N(t) = \sum_{\sigma \in \Sigma_N} \exp(-H_{N,t}(\sigma)),$$

and for any $f : \Sigma_N^n \to \mathbb{R}$, set

$$\rho_t(f) = (Z_N(t))^{-n} \sum_{\sigma^1,\ldots,\sigma^n \in \Sigma_N} f(\sigma^1,\ldots,\sigma^n) \exp\left(-\sum_{l=1}^n H_{N,t}(\sigma^l)\right).$$

The following basic relation will be essential in order to get the fluctuation result on $Z_N$:

PROPOSITION 2.2. *Assume $\beta$ satisfies the following assumption:*

(H) *$\beta$ is smaller than the constant $\beta_0$ such that*

$$16\, 2^{1/2} \beta_0 e^{16\beta_0^2} = 1.$$

*Then, for any $t \in [0,1]$, we have*

$$\mathbf{E}[\rho_t((R_{1,2} - q)^2)] \leq \frac{L}{N}.$$

PROOF. In order to show this proposition, we will use a stochastic version of the cavity method: fix $t \in [0,1]$, and for $v \in [0,t]$, set

$$-H_{N,t,v}(\sigma) = -H_{N-1,t}(\sigma)$$
$$+ \frac{\beta}{N^{1/2}} \sum_{i \leq N-1} B_{i,N}(v)\sigma_i\sigma_j + (\beta q^{1/2} X_N(v) + h)\sigma_N.$$

We will call $Z_{N,t}(v)$ [resp. $\rho_{t,v}(\cdot)$] the associated normalizing constant (resp. average with respect to the Gibbs measure). Let also $f$ be a real-valued function defined on $\Sigma_N^n$ for $N, n \geq 1$. Then $\rho_{t,v}(f)$ can be seen as a deterministic $C^2$ function of $\{B_{i,N}; i \leq N-1\}$ and $X_N$. Applying Itô's formula and Lemma 2.1 to this function, we get, for any $v \in [0, t]$,

$$\partial_v(\mathbf{E}[\rho_{t,v}(f)]) = \beta^2 \sum_{1 \leq l < l' \leq n} \mathbf{E}[\rho_{t,v}(f \sigma_N^l \sigma_N^{l'}(R_{l,l'} - q))]$$

$$- \beta^2 n \sum_{l=1}^{n} \mathbf{E}[\rho_{t,v}(f \sigma_N^l \sigma_N^{n+1}(R_{l,n+1} - q))]$$

$$+ \beta^2 \frac{n(n+1)}{2} \mathbf{E}[\rho_{t,v}(f \sigma_N^{n+1} \sigma_N^{n+2}(R_{n+1,n+2} - q))].$$

Notice also that:

1. For $v = t$, the symmetry among sites property holds true.
2. For $v = 0$, if $f^- : \Sigma_N^n \to \mathbb{R}$ depends only on the $N-1$ first coordinates, and $I$ is a subset of $\{1, \ldots, n\}$, then

$$\mathbf{E}\left[\rho_{t,0}\left(f^- \prod_{l \in I} \sigma_N^l\right)\right] = \mathbf{E}[(\tanh(\beta q^{1/2} Y + h))^{|I|}] \mathbf{E}[\rho_{t,0}(f^-)].$$

Hence, the proof of the announced result can proceed as in [13], Section 2.4. □

2.2. *Central limit theorem for the free energy.* For $t \in [0, 1]$, set $p_N(t) = \frac{1}{N} \log(Z_N(t))$. This quantity is usually called the free energy of the spin system. In this section we will get the announced central limit theorem for $p_N(t)$, with a strategy that can be briefly described as follows:

(i) We will compute the evolution of $\exp(-H_{N,t}(\sigma))$.

(ii) We will get the equation followed by the renormalized fluctuations of $Z_N(t)$, namely,

$$(8) \quad Y_N(t) = N^{1/2}\left(p_N(t) - \log(2) - \frac{\beta^2 t(1-q)^2}{4} - \mathbf{E}[\log(\cosh(\beta q^{1/2} z + h))]\right).$$

(iii) We will calculate $\mathbf{E}[\exp(\iota u Y_N(t))]$ through the application of Itô's formula, for $u \in \mathbb{R}$ and $\iota \equiv (-1)^{1/2}$, and show that this quantity tends to $e^{-u^2 \eta^2/2}$ for a positive constant $\eta^2$.

Let us turn now to the first point of this program:

LEMMA 2.3. *For any $t \in [0,1]$ and $\sigma \in \Sigma_N$, we have*

$$e^{-H_{N,t}(\sigma)} = \exp\left(\sum_{i \leq N} \sigma_i(\beta q^{1/2}\eta_i + h)\right)$$

$$+ \frac{\beta}{N^{1/2}} \sum_{i<j} \sigma_i \sigma_j \int_0^t e^{-H_{N,s}(\sigma)} dB_{i,j}(s)$$

$$+ \beta q^{1/2} \sum_{i \leq N} \sigma_i \int_0^t e^{-H_{N,s}(\sigma)} dW_i(s)$$

$$- \beta q^{1/2} \sum_{i \leq N} \sigma_i \int_0^t e^{-H_{N,s}(\sigma)} \frac{X_i(s)}{1-s} ds$$

$$+ \frac{N\beta^2}{4}\left(\frac{N-1}{N} + 2q\right) \int_0^t e^{-H_{N,s}(\sigma)} ds.$$

PROOF. The exponential being a $C^2$ function, thanks to relation (7) and Itô's formula, we get

$$e^{-H_{N,t}(\sigma)} = e^{-H_{N,0}(\sigma)} - \int_0^t e^{-H_{N,s}(\sigma)} dH_{N,s}(\sigma)$$

$$+ \tfrac{1}{2} \int_0^t e^{-H_{N,s}(\sigma)} d\langle H_{N,\cdot}(\sigma)\rangle_s.$$

In order to get the announced result, it is enough to compute the quantity $\langle H_{N,\cdot}(\sigma)\rangle_t$. But, invoking the fact that all finite variation processes have a null quadratic variation, and using the independence of all the $B_{i,j}, W_i$, we have

$$\langle H_{N,\cdot}(\sigma)\rangle_t = \frac{\beta^2}{N}\sum_{i<j}(\sigma_i\sigma_j)^2 + \beta^2 q \sum_{i \leq N}(\sigma_i)^2 = \frac{N\beta^2}{2}\left(\frac{N-1}{N} + 2q\right). \quad \square$$

We can now get the differential of $Y_N(t)$, where $Y_N$ is defined by (8). Define first, for $x \in \mathbb{R}$, the function $\Phi$ by

$$\Phi(x) = \log(\cosh(\beta q^{1/2} x + h)).$$

Then we have the following:

PROPOSITION 2.4. *The quantity $Y_N(t)$ defined by (8) satisfies, for $t \in [0,1]$,*

$$Y_N(t) = U_N + \sum_{m \leq 2} M_{m,N}(t) - (V_{1,N}(t) - V_{2,N}(t)) - V_{3,N}(t),$$

*with*

$$U_N = N^{1/2}\left(\frac{1}{N}\sum_{i\leq N} \Phi(\eta_i) - \mathbf{E}[\Phi(Y)]\right),$$

$$M_{1,N}(t) = \frac{\beta}{N}\sum_{i<j}\int_0^t \rho_s(\sigma_i\sigma_j)\,dB_{i,j}(s),$$

$$M_{2,N}(t) = \frac{\beta q^{1/2}}{N^{1/2}}\sum_{i\leq N}\int_0^t \rho_s(\sigma_i)\,dW_i(s),$$

$$V_{1,N}(t) = \frac{\beta q^{1/2}}{N^{1/2}}\sum_{i\leq N}\int_0^t \rho_s(\sigma_i)\frac{X_i(s)}{1-s}\,ds,$$

$$V_{2,N}(t) = \beta^2 q \sum_{i\leq N}\int_0^t (1-\rho_s(\sigma_i^1\sigma_i^2))\,ds,$$

$$V_{3,N}(t) = \frac{\beta^2 N^{1/2}}{4}\int_0^t \rho_s((R_{1,2}-q)^2)\,ds.$$

REMARK 2.5. As we will see in the proof of Proposition 2.4, $V_{2,N}$ is obtained as the term making $V_{1,N}(t) - V_{2,N}(t)$ a zero mean random variable.

PROOF OF PROPOSITION 2.4. Recall that $Z_N(t) = \sum_{\sigma\in\Sigma_N}\exp(-H_{N,t}(\sigma))$ is almost surely a strictly positive random variable. Hence, Itô's formula can be applied to obtain

(9)  $$\log(Z_N(t)) = \log(Z_N(0)) + A_N(t),$$

with

$$A_N(t) = \int_0^t (Z_N(s))^{-1}\,dZ_N(s) - \tfrac{1}{2}\int_0^t (Z_N(s))^{-2}\,d\langle Z_N\rangle_s.$$

It is easily checked that

$$\log(Z_N(0)) = N\log(2) + \sum_{i\leq N}\log(\cosh(\beta q^{1/2}\eta_i + h)).$$

Now, Lemma 2.3 yields that $Z_N(t)$ is a continuous semi-martingale, whose martingale part is

$$\hat{M}_N(t) = \frac{\beta}{N^{1/2}}\sum_{i<j}\sum_{\sigma\in\Sigma_N}\sigma_i\int_0^t e^{-H_{N,s}(\sigma)}\,dB_{i,j}(s)$$

$$+ \beta q^{1/2}\sum_{i\leq N}\sum_{\sigma\in\Sigma_N}\sigma_i\int_0^t e^{-H_{N,s}(\sigma)}\,dW_i(s).$$

Hence,
$$\int_0^t (Z_N(s))^{-2} d\langle Z_N \rangle_s$$
$$= \int_0^t (Z_N(s))^{-2} d\langle \hat{M}_N \rangle_s$$
$$= \frac{\beta^2}{N} \sum_{i<j} \int_0^t \rho_s(\sigma_i^1 \sigma_j^1 \sigma_i^2 \sigma_j^2) \, ds + \beta^2 q \sum_{i \leq N} \int_0^t \rho_s(\sigma_i^1 \sigma_i^2) \, ds.$$

Observing that

(10)
$$\sum_{i<j} \rho_s(\sigma_i^1 \sigma_j^1 \sigma_i^2 \sigma_j^2) = \frac{N^2}{2} \left( \rho_s(R_{1,2}^2) - \frac{1}{N} \right),$$
$$\sum_{i \leq N} \rho_s(\sigma_i^1 \sigma_i^2) = N \rho_s(R_{1,2})$$

yields the following identity:
$$\int_0^t (Z_N(s))^{-2} d\langle Z_N \rangle_s = \frac{\beta^2 N}{2} \int_0^t \left( \rho_s(R_{1,2}^2) + 2q\rho_s(R_{1,2}) - \frac{1}{N} \right) ds.$$

The same kind of considerations can be used to obtain an expression for $(Z_N(s))^{-1} dZ_N(s)$. Then, plugging those relations into (9), using Lemma 2.3 and a little algebra, we get

$$A_N(t) = \frac{\beta}{N^{1/2}} \sum_{i<j} \int_0^t \rho_s(\sigma_i \sigma_j) \, dB_{i,j}(s)$$
$$+ \beta q^{1/2} \sum_{i \leq N} \int_0^t \rho_s(\sigma_i) \, dW_i(s)$$
$$- \beta q^{1/2} \sum_{i \leq N} \int_0^t \rho_s(\sigma_i) \frac{X_i(s)}{1-s} \, ds$$
$$- \frac{\beta^2 N}{4} \int_0^t [(\rho_s(R_{1,2}^2) - 1) + 2q(\rho_s(R_{1,2}) - 1)] \, ds.$$

Let us center now on the third term of this sum: according to (6), we have, for any $s \in [0,1]$, with obvious notation,
$$\mathbf{E}\left[ \rho_s(\sigma_i) \frac{X_i(s)}{1-s} \right] = \mathbf{E}[\partial_{X_i(s)}(\rho_s(\sigma_i))].$$

But
$$\rho_s(\sigma_i) = \frac{\sum_{\sigma^1 \in \Sigma_N} \sigma_i^1 \exp(-\hat{H}_{N,s}(\sigma^1) + \beta q^{1/2} \sum_{i \leq N} X_i(s) \sigma_i^1)}{\sum_{\sigma^2 \in \Sigma_N} \exp(-\hat{H}_{N,s}(\sigma^2) + \beta q^{1/2} \sum_{i \leq N} X_i(s) \sigma_i^2)},$$

where
$$-\hat{H}_{N,s}(\sigma) = \frac{\beta}{N^{1/2}} \sum_{i<j} B_{i,j}(s)\sigma_i\sigma_j + h \sum_{i \leq N} \sigma_i,$$

and it is easily seen from that expression that

(11) $$\partial_{X_i(s)}(\rho_s(\sigma_i)) = \beta q^{1/2} \rho_s(1 - \sigma_i^1 \sigma_i^2).$$

Hence,

$$\log(Z_N(t)) = N\log(2) + \sum_{i \leq N} \Phi(\eta_i) + \frac{N\beta^2(1-q)^2 t}{4}$$
$$+ N^{1/2}[M_{1,N}(t) + M_{2,N}(t) - (V_{1,N}(t) - V_{2,N}(t)) + V_{3,N}(t)].$$

Now, we get the desired result by substracting

$$\log(2) + \frac{N\beta^2(1-q)^2 t}{4} + \mathbf{E}[\Phi(Y)],$$

and renormalizing this expression. □

We are now ready to state the main result of this section:

THEOREM 2.6. *Let $Y_N(t)$ be defined by (8) for $t \in [0, 1]$. Then $Y_N(1)$ converges in distribution to a $\mathcal{N}(0, \tau^2)$ variable, with*

$$\tau^2 = \nu^2 - \frac{\beta^2 q^2}{2} \qquad \text{where } \nu^2 = \text{Var}(\log(\cosh(\beta q^{1/2} Y + h))).$$

PROOF. Let $u \in \mathbb{R}$. Once we know the differential of $X_N$, given by Proposition 2.4, we can apply Itô's formula to the (complex-valued) function $x \mapsto e^{\iota u x}$ to obtain, for $t \in [0, 1]$,

$$e^{\iota u Y_N(t)} = D_1(N) + \sum_{m=2}^{6} D_{m,N}(t),$$

with

$$D_1(N) = e^{\iota u U_N},$$

$$D_{2,N}(t) = \frac{\iota u \beta}{N} \sum_{i<j} \int_0^t e^{\iota u Y_N(s)} \rho_s(\sigma_i \sigma_j) \, dB_{i,j}(s)$$
$$+ \frac{\iota u \beta q^{1/2}}{N^{1/2}} \sum_{i \leq N} \int_0^t e^{\iota u Y_N(s)} \rho_s(\sigma_i) \, dW_i(s),$$

$$D_{3,N}(t) = -\frac{\iota u \beta q^{1/2}}{N^{1/2}} \sum_{i \leq N} \int_0^t e^{\iota u Y_N(s)} \rho_s(\sigma_i) \frac{X_i(s)}{1-s} \, ds$$

and

$$D_{4,N}(t) = \frac{\iota u \beta^2 q}{N^{1/2}} \sum_{i \leq N} \int_0^t e^{\iota u Y_N(s)} [1 - \rho_s(\sigma_i^1 \sigma_i^2)] \, ds,$$

$$D_{5,N}(t) = -\frac{\iota u \beta^2 N^{1/2}}{4} \int_0^t e^{\iota u Y_N(s)} \rho_s((R_{1,2} - q)^2) \, ds,$$

$$D_{6,N}(t) = -\frac{u^2 \beta^2}{2N^2} \sum_{i<j} \int_0^t e^{\iota u Y_N(s)} \rho_s(\sigma_i^1 \sigma_j^1 \sigma_i^2 \sigma_j^2) \, ds,$$

$$D_{7,N}(t) = -\frac{u^2 \beta^2 q}{2} \int_0^t e^{\iota u Y_N(s)} \rho_s(R_{1,2}) \, ds.$$

Denote also $\mathbf{E}[e^{\iota u Y_N(t)}]$ by $\psi_{N,u}(t)$. We will now estimate the previous terms separately (in the sequel, $\kappa$ will stand for a positive constant that can change from line to line): first, from the classical central limit theorem for i.i.d. random variables, we have

(12) $\quad D_1(N) = e^{-\nu^2 u^2/2} + \mathcal{R}_1(N), \qquad |\mathcal{R}_1(N)| \leq \dfrac{\kappa}{N^{1/2}}.$

The term $D_{2,N}(t)$ being of zero mean, we will turn now to the estimation of $\mathbf{E}[D_{3,N}(t)]$: just as in the proof of Proposition 2.4, an integration by parts yields

$$\mathbf{E}[D_{3,N}(t)] = -\frac{\iota u \beta q^{1/2}}{N^{1/2}} \sum_{i \leq N} \int_0^t \mathbf{E}[\partial_{X_i(s)} (e^{\iota u Y_N(s)} \rho_s(\sigma_i))] \, ds$$

$$= \frac{u^2 \beta q^{1/2}}{N^{1/2}} \sum_{i \leq N} \int_0^t \mathbf{E}[(\partial_{X_i(s)} Y_N(s)) e^{\iota u Y_N(s)} \rho_s(\sigma_i)] \, ds - \mathbf{E}[D_{4,N}(t)].$$

Observe now that, following the lines of the proof of (11), we have

$$\partial_{X_i(s)} Y_N(s) = \frac{\partial_{X_i(s)} Z_N(s)}{N^{1/2} Z_N(s)} = \frac{\beta q^{1/2}}{N^{1/2}} \rho_s(\sigma_i),$$

and, hence,

$$\mathbf{E}[D_{3,N}(t) + D_{4,N}(t)] = u^2 \beta^2 q \int_0^t \mathbf{E}[e^{\iota u Y_N(s)} \rho_s(R_{1,2})] \, ds.$$

It is now easily seen, thanks to Proposition 2.2, that

(13) $\quad \mathbf{E}[D_{3,N}(t) + D_{4,N}(t)] = (u \beta q)^2 \int_0^t \psi_{N,u}(s) \, ds + \mathcal{R}_{3,t}(N),$

with

$$|\mathcal{R}_{3,t}(N)| \leq \frac{\kappa}{N^{1/2}},$$

uniformly in $t \in [0, 1]$. A direct application of Proposition 2.2 also gives that, for any $t \in [0, 1]$,

$$\mathbf{E}[D_{5,N}(t)] \leq \frac{\kappa}{N^{1/2}}. \tag{14}$$

Invoking now (10), we have

$$\mathbf{E}[D_{6,N}(t)] = \frac{u^2 \beta^2}{4} \int_0^t e^{\iota u Y_N(s)} \left( \rho_s(R_{1,2}^2) - \frac{1}{N} \right) ds,$$

and by Proposition 2.2 again, we get

$$\mathbf{E}[D_{6,N}(t)] = -\frac{(u\beta q)^2}{4} \int_0^t \psi_{N,u}(s) \, ds + \mathcal{R}_{6,t}(N),$$

$$|\mathcal{R}_{6,t}(N)| \leq \frac{\kappa}{N^{1/2}}, \tag{15}$$

and the same type of arguments give

$$\mathbf{E}[D_{7,N}(t)] = -\frac{(u\beta q)^2}{2} \int_0^t \psi_{N,u}(s) \, ds + \mathcal{R}_{7,t}(N),$$

$$|\mathcal{R}_{7,t}(N)| \leq \frac{\kappa}{N^{1/2}}. \tag{16}$$

Putting together (12)–(16), we have, finally,

$$\psi_{N,u}(t) = e^{-\nu^2 u^2/2} + \frac{(u\beta q)^2}{4} \int_0^t \psi_{N,u}(s) \, ds + \hat{\mathcal{R}}_{N,u}(t),$$

with

$$|\hat{\mathcal{R}}_{N,u}(t)| \leq \frac{\kappa}{N^{1/2}},$$

which ends the proof by a Gronwall type argument. □

**3. The perceptron case.** The aim of this section is to show that our method can be applied in various contexts. We will illustrate this point by showing a CLT for the free energy of another canonical model of spin glasses, namely, the perceptron model.

3.1. *Statement of the results.* The perceptron model is still based on the state space $\Sigma_N = \{-1, 1\}^N$, for $N \geq 1$, and can be described as follows: consider a positive integer $M = \alpha N$ such that $M = \alpha N$ for a given $\alpha > 0$; $u$ will be a continuous function defined on $\mathbb{R}$ satisfying $|u| \leq D$ for a strictly positive constant $D$ (some additional assumption will be made on $u$ further on), and $\{g_{i,k}; i \geq 1, k \geq 1\}$ is a family of independent standard Gaussian random variables. The random measure we will consider on $\Sigma_N$ will be of the form

$$(Z_{N,M})^{-1}\exp(-H_{N,M}(\sigma))\mu_N(d\sigma),$$

where $\mu_N$ is the uniform measure on $\Sigma_N$, and for any $m \geq 1$,

$$-H_{N,m}(\sigma) = \sum_{k \leq m} u(S_k),$$

where, for $l \geq 1$ and a replica $\sigma^l$, we have

(17) $$S_k^l = S_k^l(\sigma^l) = \frac{1}{N^{1/2}} \sum_{i \leq N} g_{i,k}\sigma_i^l.$$

The normalization constant associated to this model is $Z_{N,M}$, with $Z_{N,m}$ defined, for $m \geq 1$, by

$$Z_{N,m} = \sum_{\sigma \in \Sigma_N} \exp(-H_{N,m}(\sigma)),$$

and for any $M \geq 1$, $f:\Sigma_N^n \to \mathbb{R}$, $m \leq M$, the average of $f$ with respect to the measure defined by $H_{N,m}$ will be

$$\rho_m(f) = Z_{N,m}^{-n} \sum_{\sigma^1,\ldots,\sigma^n \in \Sigma_N} f(\sigma^1,\ldots,\sigma^n) \exp\left(-\sum_{l=1}^n H_{N,m}(\sigma^l)\right).$$

The analysis of the perceptron model relies, as in the SK case, on the limiting behavior of the overlap $R_{1,2}$, and we recall the following result, taken from [13], Chapter 3: for $m \leq M$ and $\alpha_m = \frac{m}{N}$, consider the system of equations

(18) $$q = \mathbf{E}[\tanh^2(r^{1/2}z)], \qquad r = \alpha_m \mathbf{E}[\Psi^2(q^{1/2}z,(1-q)^{1/2})],$$

where $z, \xi$ are two independent standard normal random variables, and $\Psi$ is defined by

$$\Psi(x,y) = \frac{\mathbf{E}[\xi \exp(u(x+\xi y))]}{y\mathbf{E}[\exp(u(x+\xi y))]}.$$

Then we have the following:

PROPOSITION 3.1. *Assume $u$ and $\alpha$ satisfy:*

(H1) $\|u\| \leq D$ *and, for $L > 0$ large enough, $L\alpha \exp(LD) \leq 1$.*
(H2) *There exists a positive constant $L^*$ and a small enough constant $c_3$ such that, for any $l \leq 6$, $|u^{(l)}| \leq L^* e^{L^*D} e^{c_3 N}$.*

*Then, for any $m \leq M$,*

1. The system (18) has a unique solution $(q_m, r_m) \in [0,1]^2$.
2. The following $L^2$ convergence for $R_{1,2}$ holds true:

$$\mathbf{E}[\rho_m((R_{1,2} - q_m)^2)] \leq \frac{\kappa}{N}.$$

3. Set $p_{N,m} = \frac{1}{N}\mathbf{E}[\log(Z_{N,m})]$ and

$$\Phi(m) = \log(2) + \mathbf{E}[\log(\cosh(r_m^{1/2} Y))] - \tfrac{1}{2} r_m (1 - q_m)$$
$$+ \alpha_m \mathbf{E}[\log(\hat{\mathbf{E}}[\exp(u(q_m^{1/2}\eta + (1-q_m)^{1/2}\hat{\eta}))])],$$

where $Y, \eta, \hat{\eta}$ are independent $\mathcal{N}(0,1)$ random variables, and $\hat{\mathbf{E}}$ designates the expectation with respect to $\hat{\eta}$. Then

$$|p_{N,m} - \Phi(m)| \leq \frac{\kappa}{N}.$$

In the sequel of the section, we will use the following auxiliary random variables, for $l \geq 1$ and $1 \leq m \leq M$:

$$\theta_m^l = q_m^{1/2} \eta + (1 - q_m)^{1/2} \hat{\eta}^l,$$

where $\eta$ and $\{\hat{\eta}^l; l \geq 1\}$ are independent standard Gaussian random variables. We will denote by $\hat{\mathbf{E}}$ the expectation with respect to the random variables $\hat{\eta}^l$, and also set $\theta_m = \theta_m^1$ when only one replica is considered. Associated to those random variables $\theta_m$ will be the quantity $\xi_m$, defined by

(19) $$\xi_m = \log(\hat{\mathbf{E}}[\exp(u(\theta_m))]).$$

It is then easily checked that

(20) $$\operatorname{Var}(\xi_m) = Q\left(\frac{m}{N}\right), \qquad \sup_{m \leq M} \mathbf{E}[\xi_m^4] < \infty,$$

where $Q : [0, \alpha] \to \mathbb{R}_+$ is a $C^1$ function.

Our aim, in this part of the paper, is to prove the following theorem by means of the stochastic calculus method:

THEOREM 3.2. *Under assumptions* (H1) *and* (H2), *as $N$ tends to infinity,*

$$N^{1/2}\left[\frac{\log(Z_{N,M})}{N} - \Phi(M)\right] \xrightarrow{(\mathcal{L})} Y,$$

*where $Y$ is a Gaussian random variable with variance*

$$\tau^2 = \frac{1}{\alpha}\int_0^\alpha Q(x)\, dx,$$

*with $Q$ defined at* (20).

3.2. *Proof of the CLT.* Let us begin the proof of our result by the following elementary property, that we label for further use:

LEMMA 3.3. *Let $1 \le m \le M$, and set*
$$\Delta_{m-1,m}\Phi = \Phi(m) - \Phi(m-1).$$
*Then*
$$\Delta_{m-1,m}\Phi - \frac{\mathbf{E}[\xi_m]}{N} = \mathcal{R}_{m,N},$$
*with $\mathcal{R}_{m,N} \le \frac{\kappa}{N^2}$ for a positive constant $\kappa$, and where $\xi_m$ is defined by* (19).

PROOF. The lector is referred to [13], Theorem 3.4.2, for a complete proof of this fact. Let us just mention that this mysterious relation relies on the fact that, since $(q_m, r_m)$ is the solution to (18), then $\Phi(m)$ is a function of $\alpha = \frac{m}{N}$, say $\hat{\Phi}(\alpha)$, and it can be shown by elementary computations that
$$\partial_\alpha \hat{\Phi}(\alpha) = \mathbf{E}[\log(\hat{\mathbf{E}}[e^{u(\theta)}])]. \qquad \square$$

We will now write the quantity $\frac{1}{N}\log(Z_{N,M}) - \Phi(M)$ in a convenient way in order to compute its Fourier transform: it is easily checked that, for $m \le M$,
$$\frac{1}{N}\log(Z_{N,M}) - \Phi(M) = \frac{1}{N}\log(Z_{N,M-1}) - \Phi(M-1)$$
$$+ (\log(\rho_{M-1}(e^{u(S_M)})) - \Delta_{M-1,M}\Phi),$$
and iterating this decomposition, we get
$$(21) \qquad \frac{1}{N}\log(Z_{N,M}) - \Phi(M) = \sum_{m \le M} Y_m,$$
where
$$Y_m = \log(\rho_{m-1}(e^{u(S_m)})) - \Delta_{m-1,m}\Phi.$$
Our proof of the CLT will be based on getting some information on the characteristic function of each $Y_m$ separately. This will be achieved by changing the random variable $S_m^l$ along the stochastic path
$$S_m^l(t) = \frac{1}{N^{1/2}} \sum_{i \le N} B_{i,m}(t)\sigma_i^l + q_m^{1/2} X_m(t) + (1-q_m)^{1/2} \hat{X}_m^l(t),$$
defined for $t \in [0,1]$, where $X_m, \hat{X}_m^l$ are reversed time Brownian motions, solutions to
$$X_m(t) = \eta_m - \int_0^t \frac{X_m(s)}{1-s} ds + W_m(t),$$
$$\hat{X}_m^l = \hat{\eta}_m^l - \int_0^t \frac{\hat{X}_m^l(s)}{1-s} ds + \hat{W}_m^l(t),$$

for $t \in [0,1]$, with some independent standard Brownian motions

$$\{B_{i,m}, i, m \geq 1\}, \qquad \{W_m, m \geq 1\}, \qquad \{\hat{W}_m^l, m, l \geq 1\}.$$

We will also denote by $\mathbf{E}_m$ the expectation conditioned on all the randomness in the $g_{i,j}, B_{i,j}, X_j, \ldots$, for $j \leq m$, and we still call $\hat{\mathbf{E}}$ the expectation with respect to the random variables of the form $\hat{\eta}, \hat{X}, \hat{W}$. Next, we set

(22) $$Y_m(t) = \frac{1}{N} \log(\hat{\mathbf{E}}[\rho_{m-1}(e^{u(S_m(t))})]) - \Delta_{m-1,m}\Phi.$$

Eventually, define $\bar{\xi}_m = \xi_m - \mathbf{E}[\xi_m]$. Notice that all those quantities are related by

$$S_m^l = S_m^l(1), \qquad \theta_m^l = S_m^l(0), \qquad \xi_m = \log(\hat{\mathbf{E}}[\rho_{m-1}(e^{u(S_m(0))})]).$$

With all those notations in mind, we can now gather some information about the Fourier transform of $Y_m$.

PROPOSITION 3.4. *For any $m \leq M$ and $v \in \mathbb{R}$, we have*

$$\mathbf{E}[e^{\iota v N^{1/2} \sum_{k \leq m} Y_k}] = \mathbf{E}[e^{\iota v N^{1/2} \bar{\xi}_m}]\mathbf{E}[e^{\iota v N^{1/2} \sum_{k \leq m-1} Y_k}] + \mathcal{R}_{m,N},$$

*with* $|\mathcal{R}_{m,N}| \leq \frac{\kappa}{N^{3/2}}$.

PROOF. Notice first that

(23) $$\mathbf{E}[e^{\iota v N^{1/2} \sum_{k \leq m} Y_k}] = \mathbf{E}[e^{\iota v N^{1/2} \sum_{k \leq m-1} Y_k} \mathbf{E}_{m-1}[e^{\iota v N^{1/2} Y_m}]].$$

We will try now to get an expansion of $e^{\iota v N^{1/2} Y_m}$ using Itô's formula.

STEP 1. With $Y_m(t)$ defined by (22), recall that we have $Y_m = Y_m(1)$ and

$$\mathbf{E}_{m-1}[e^{\iota v N^{1/2} Y_m(0)}] = \mathbf{E}_{m-1}[e^{\iota v N^{1/2}(\xi_m - \Delta_{m-1,m}\Phi)}].$$

Hence, on account of Lemma 3.3, we get

$$|\mathbf{E}_{m-1}[e^{\iota v N^{1/2} \bar{\xi}_m}] - \mathbf{E}_{m-1}[e^{\iota v N^{1/2} Y_m(0)}]|$$
$$\leq |\mathbf{E}_{m-1}[e^{\iota v N^{1/2} \bar{\xi}_m}(1 - \exp(\iota v N^{1/2}(\Delta_{m-1,m}\Phi - \mathbf{E}[\xi_m])))]| \leq \frac{\kappa}{N^{3/2}}.$$

STEP 2. The evolution in time of $Y_m(t)$ can be described as follows: first, setting

$$U_1 = u'e^u, \qquad U_2 = ((u')^2 + u'')e^u,$$

we get, by a simple application of Itô's formula,

$$e^{u(S_m(t))} = e^{u(\theta)} + \frac{1}{N^{1/2}} \sum_{i \leq N} \sigma_i \int_0^t U_1(S_m(s)) \, dB_{i,m}(s)$$

$$+ q_m^{1/2} \int_0^t U_1(S_m(s)) \left[ dW_m(s) - \frac{X_m(s)}{1-s} \, ds \right]$$

$$+ (1 - q_m)^{1/2} \int_0^t U_1(S_m(s)) \left[ d\hat{W}_m(s) - \frac{\hat{X}_m(s)}{1-s} \, ds \right]$$

$$+ \int_0^t U_2(S_m(s)) \, ds,$$

which gives immediately

$$\rho_{m-1}(e^{u(S_m(t))}) = e^{u(\theta)} + \frac{1}{N^{1/2}} \sum_{i \leq N} \int_0^t \rho_{m-1}(\sigma_i U_1(S_m(s))) \, dB_{i,m}(s)$$

$$+ q_m^{1/2} \int_0^t \rho_{m-1}(U_1(S_m(s))) \left[ dW_m(s) - \frac{X_m(s)}{1-s} \, ds \right]$$

$$+ (1 - q_m)^{1/2} \int_0^t \rho_{m-1}(U_1(S_m(s))) \left[ d\hat{W}_m(s) - \frac{\hat{X}_m(s)}{1-s} \, ds \right]$$

$$+ \int_0^t \rho_{m-1}(U_2(S_m(s))) \, ds.$$

Let us now compute the expectation with respect to the randomness in $\hat{X}_m$: by integration by parts with respect to $\frac{\hat{X}_m(s)}{1-s}$, we get

$$\hat{\mathbf{E}}[\rho_{m-1}(e^{u(S_m(t))})] = \hat{\mathbf{E}}[e^{u(\theta)}] + \frac{1}{N^{1/2}} \sum_{i \leq N} \int_0^t \hat{\mathbf{E}}[\rho_{m-1}(\sigma_i U_1(S_m(s)))] \, dB_{i,m}(s)$$

$$+ q_m^{1/2} \int_0^t \hat{\mathbf{E}}[\rho_{m-1}(U_1(S_m(s)))] \left[ dW_m(s) - \frac{X_m(s)}{1-s} \, ds \right]$$

$$+ q_m \int_0^t \hat{\mathbf{E}}[\rho_{m-1}(U_2(S_m(s)))] \, ds.$$

Now, like in the SK case, we will apply the function $f(x) = e^{\iota v N^{1/2} \log(x)}$ to the quantity

$$Z_{m-1}(t) \equiv \hat{\mathbf{E}}[\rho_{m-1}(e^{u(S_m(t))})],$$

and, according to Itô's formula, we get

$$e^{\iota v N^{1/2} Y_m(t)} = e^{\iota v N^{1/2} (\xi_m - \Delta_{m-1,m} \Phi)} + \sum_{j \leq 3} A_j(t) + \sum_{j \leq 2} B_j(t),$$

with

$$A_1(t) = \frac{\iota v}{N} \sum_{i \leq N} \int_0^t e^{\iota v N^{1/2} Y_m(s)} \frac{\hat{\mathbf{E}}[\rho_{m-1}(\sigma_i U_1(S_m(s)))]}{Z_{m-1}(s)} dB_{i,m}(s),$$

$$A_2(t) = \frac{\iota v q_m^{1/2}}{N^{1/2}} \int_0^t e^{\iota v N^{1/2} Y_m(s)} \frac{\hat{\mathbf{E}}[\rho_{m-1}(U_1(S_m(s)))]}{Z_{m-1}(s)} \left[ dW_m(s) - \frac{X_m(s)}{1-s} ds \right],$$

$$A_3(t) = \frac{\iota v q_m}{N^{1/2}} \int_0^t e^{\iota v N^{1/2} Y_m(s)} \frac{\hat{\mathbf{E}}[\rho_{m-1}(U_2(S_m(s)))]}{Z_{m-1}(s)} ds,$$

and where the quantities $B_1(t)$ and $B_2(t)$ are given by

$$B_1(t) = \frac{\iota v}{2N^{1/2}} \int_0^t e^{\iota v N^{1/2} Y_m(s)} \frac{\hat{\mathbf{E}}[\rho_{m-1}((R_{1,2}+q)U_1(S_m^1(s))U_1(S_m^2(s)))]}{(Z_{m-1}(s))^2} ds$$

and

$$B_2(t) = \frac{\iota v^2}{2N} \int_0^t e^{\iota v N^{1/2} Y_m(s)} \frac{\hat{\mathbf{E}}[\rho_{m-1}((R_{1,2}-q)U_1(S_m^1(s))U_1(S_m^2(s)))]}{(Z_{m-1}(s))^2} ds.$$

Taking now the conditional expectation $\mathbf{E}_{m-1}$, and integrating by parts with respect to $\frac{X_m(s)}{1-s}$ [notice that one has to differentiate, $e^{\iota v N^{1/2} Y_m(s)}$, $U_1(S_m(s))$ and $Z_{m-1}(s)$], the above expression simplifies into

(24) $$\mathbf{E}_{m-1}[e^{\iota v N^{1/2} Y_m(t)}] = I_1 - I_2(t) - I_3(t),$$

with

$$I_1(t) = e^{\iota v N^{1/2}(\xi_m - \Delta_{m-1,m}\Phi)},$$

the quantities $I_2(t)$ and $I_3(t)$ being defined by

$$I_2(t) = \frac{\iota v q_m}{N^{1/2}} \int_0^t \mathbf{E}_{m-1} \left[ e^{\iota v N^{1/2} Y_m(s)} \frac{\rho_{m-1}((R_{1,2}-q)U_1(S_m^1(s))U_1(S_m^2(s)))}{(Z_{m-1}(s))^2} \right] ds$$

and

$$I_3(t) = \frac{\iota v^2}{2N} \int_0^t \mathbf{E}_{m-1} \left[ e^{\iota v N^{1/2} Y_m(s)} \frac{\rho_{m-1}((R_{1,2}-q)U_1(S_m^1(s))U_1(S_m^2(s)))}{(Z_{m-1}(s))^2} \right] ds.$$

STEP 3. Let us go back to expression (23), and note that, for any $t \in [0,1]$, invoking relation (24), we have

(25) $$\mathbf{E}[e^{\iota v N^{1/2}(\sum_{k \leq m-1} Y_k + Y_m(t))}] = \mathbf{E}[e^{\iota v N^{1/2} \sum_{k \leq m-1} Y_k}(I_1 - I_2(t) - I_3(t))].$$

Let us analyze now the three terms we have obtained: we have already shown (see Step 1) that

(26) $$|I_1 - \mathbf{E}[e^{\iota v N^{1/2} \bar{\xi}_m}]| \leq \frac{\kappa}{N^{3/2}},$$

which easily yields

$$|\mathbf{E}[e^{\iota v N^{1/2}\sum_{k\leq m-1} Y_k} I_1] - \mathbf{E}[e^{\iota v N^{1/2}\bar{\xi}_m}]\mathbf{E}[e^{\iota v N^{1/2}\sum_{k\leq m-1} Y_k}]| \leq \frac{\kappa}{N^{3/2}}.$$

On the other hand,

$$(27) \quad |\mathbf{E}[e^{\iota v N^{1/2}\sum_{k\leq m-1} Y_k} I_3(t)]| = \frac{v^2}{2N}\left|\int_0^t \mathbf{E}[e^{\iota v N^{1/2}(\sum_{k\leq m-1} Y_k + Y_m(s))} K_m(s)]\,ds\right|,$$

with

$$K_m(s) = \frac{\rho_{m-1}((R_{1,2} - q)U_1(S_m^1(s))U_1(S_m^2(s)))}{(Z_{m-1}(s))^2},$$

and the right-hand side of (27) is bounded by

$$\frac{v^2}{2N}\int_0^t \mathbf{E}^{1/2}[(R_{1,2} - q)^2]\mathbf{E}^{1/2}\left[(Z_{m-1}(s))^{-4}\prod_{j\leq 4} U_1(S_m^j(s))\right] ds.$$

Observe that the last of those terms seems to be potentially huge, since the derivatives of $u$ are allowed to grow at exponential speed with $N$. However, under assumptions (H1) and (H2), it is shown in [13], Lemma 3.3.4, by some integration by parts arguments, that this kind of term is uniformly bounded in $N$. Hence, invoking the $L^2$ estimate for $R_{1,2}$ given at Proposition 3.1, for any $t \in [0,1]$,

$$(28) \qquad |\mathbf{E}[e^{\iota v N^{1/2}\sum_{k\leq m-1} Y_k} I_3(t)]| \leq \frac{\kappa}{N^{3/2}}.$$

We still have to handle the term

$$\mathbf{E}[e^{\iota v N^{1/2}\sum_{k\leq m-1} Y_k} I_2(t)].$$

This is achieved by a long and tedious variation of Theorem 3.5.2 and Proposition 3.5.1 in [13]. Since no new idea is required, we will only sketch the proof of that step: first, one has to show that

$$(29) \quad \begin{aligned} &\mathbf{E}[e^{\iota v N^{1/2}\sum_{k\leq m-1} Y_k} I_2(t)] \\ &= \frac{\iota v q_m \gamma_m t}{N^{1/2}}\mathbf{E}[e^{\iota v N^{1/2}\sum_{k\leq m-1} Y_k} \rho_{m-1}(R_{1,2} - q)] + \mathcal{R}_{m,N}, \end{aligned}$$

with

$$\gamma_m = \mathbf{E}\left[e^{\iota v N^{1/2}\bar{\xi}_m}\frac{U_1(\theta^1)U_1(\theta^2)}{(\hat{\mathbf{E}}[e^{u(\theta)}])^2}\right], \qquad |\mathcal{R}_{m,N}| \leq \frac{\kappa}{N^{3/2}}.$$

This is obtained by a variation of the proof of [13], Theorem 3.5.2. The main differences with this last proof are that:

1. One has to use the symmetry among sites to express $\rho_{m-1}(R_{1,2} - q)$ as a function of the last spin only. Then it can be realized that the equivalent contribution of the terms II, III and IV in [13], Theorem 3.5.2, are of order $N^{-1}$ when the function $f$ considered there is a function of the last spin.
2. We also have to keep track, in the computation of the derivatives, of the terms $e^{\iota v N^{1/2} Y_m(s)}$, which increases the size of the computations. However, in the end, those additional terms are all of order $N^{-1}$. Notice that all the calculations can be done again using Itô's formula.

Then, it should also be checked that

$$(30) \qquad |\mathbf{E}[e^{\iota v N^{1/2} \sum_{k \leq m-1} Y_k} \rho_{m-1}(R_{1,2} - q)]| \leq \frac{\kappa}{N}.$$

To this purpose, the proof of Proposition 3.5.1 in [13] has to be followed, keeping track again of the terms due to $e^{\iota v N^{1/2} \sum_{k \leq m-1} Y_k}$ that are still of order $N^{-1}$.

Now, (29) and (30) easily yield

$$(31) \qquad \mathbf{E}[e^{\iota v N^{1/2} \sum_{k \leq m-1} Y_k} I_2(t)] \leq \frac{\kappa}{N^{3/2}},$$

and plugging (26), (28) and (31) into (25), we get, for any $t \in [0, 1]$,

$$|\mathbf{E}[e^{\iota v N^{1/2}(\sum_{k \leq m-1} Y_k + Y_m(t))}] - \mathbf{E}[e^{\iota v N^{1/2} \bar{\xi}_m}]\mathbf{E}[e^{\iota v N^{1/2} \sum_{k \leq m-1} Y_k}]| \leq \frac{\kappa}{N^{3/2}},$$

which implies the claim of our proposition. □

Let us go back now to the main aim of this section:

PROOF OF THEOREM 3.2. Let $v$ be an arbitrary real number. We have seen [relation (21)] that

$$\frac{1}{N} \log(Z_{N,M}) - \Phi(M) = \sum_{m \leq M} Y_m.$$

Now, Proposition 3.4 gives

$$\mathbf{E}[e^{\iota v N^{1/2} \sum_{k \leq M} Y_k}] = \mathbf{E}[e^{\iota v N^{1/2} \bar{\xi}_M}]\mathbf{E}[e^{\iota v N^{1/2} \sum_{k \leq M-1} Y_k}] + \mathcal{R}_{M,N},$$

and, iterating this decomposition, we get

$$\mathbf{E}[e^{\iota v N^{1/2} \sum_{m \leq M} Y_k}] = \mathbf{E}[e^{\iota v N^{1/2} \sum_{m \leq M} \bar{\xi}_m}] + \sum_{m \leq M} \mathcal{R}_{m,N},$$

with $|\mathcal{R}_{m,N}| \leq \frac{\kappa}{N^{3/2}}$ for any $m \leq M$. Thus,

$$|\mathbf{E}[e^{\iota v N^{1/2} \sum_{m \leq M} Y_k}] - \mathbf{E}[e^{\iota v N^{1/2} \sum_{m \leq M} \bar{\xi}_m}]| \leq \frac{\kappa}{N^{1/2}},$$

for a positive constant $\kappa$. Now, the random variables $\{\bar{\xi}_m; m \leq M\}$ are centered and independent, and we have seen that

$$\mathrm{Var}(\bar{\xi}_m) = \mathrm{Var}(\xi_m) = Q\left(\frac{m}{N}\right),$$

where $Q:[0,\alpha] \to \mathbb{R}_+$ is a $C^1$ function. Hence, by Riemann sums approximation,

$$\left|\frac{1}{N}\sum_{m \leq M} \mathrm{Var}(\bar{\xi}_m) - \frac{1}{\alpha}\int_0^\alpha Q(x)\,dx\right| \leq \frac{\kappa\|Q'\|_{[0,\alpha]}}{N}.$$

Taking into account inequality (20), the end of the proof easily follows now by some classical CLT arguments for independent random variables (see, e.g., the Lindeberg–Feller theorem in [4]). $\square$

**Acknowledgments.** This paper has partially been written with the help of a fellowship grant at the *Centre de Recerca Matemàtica* (Bellaterra, Catalonia). I would like to thank this institution for its warm hospitality, as well as for the nice atmosphere it creates for good mathematical work. I would also like to thank Francis Comets, M'hammed Eddhabi, Arnaud Guillin and Michel Talagrand for some valuable conversations concerning the paper.

INSTITUT ÉLIE CARTAN
UNIVERSITÉ HENRI POINCARÉ (NANCY)
BP 239
54506-VANDOEUVRE-LÈS-NANCY
FRANCE
E-MAIL: tindel@iecn.u-nancy.fr